\documentclass{amsart}
\usepackage{amssymb, amsmath, mathrsfs, verbatim, url, bbm}

\theoremstyle{plain}
\newtheorem{theorem}{Theorem}[section]
\newtheorem{lemma}[theorem]{Lemma}
\newtheorem{proposition}[theorem]{Proposition}
\newtheorem{corollary}[theorem]{Corollary}

\theoremstyle{definition}
\newtheorem{definition}[theorem]{Definition}
\newtheorem{example}[theorem]{Example}

\theoremstyle{remark}

\newtheorem{question}[theorem]{Question}

\numberwithin{equation}{section}

\DeclareMathOperator{\dom}{dom}
\DeclareMathOperator{\ran}{ran}
\DeclareMathOperator{\cf}{cf}

\newcommand{\nbd}{\nobreakdash}

\newcommand{\powset}[1]{\mathcal{P}(#1)}
\newcommand{\card}[1]{\lvert #1\rvert}

\newcommand{\la}{\langle \nolinebreak[4]}
\newcommand{\ra}{\nolinebreak[4] \rangle}

\newcommand{\cardc}{\mathfrak{c}}
\newcommand{\cardp}{\mathfrak{p}}
\newcommand{\cardu}{\mathfrak{u}}
\newcommand{\inv}[1]{#1^{-1}}

\newcommand{\restrict}{\upharpoonright}

\newcommand{\wma}{we may assume}
\newcommand{\Wma}{We may assume}

\newcommand{\elemsub}{\prec}

\newcommand{\mcA}{\mathcal{ A}}
\newcommand{\mcB}{\mathcal{ B}}
\newcommand{\mcC}{\mathcal{ C}}
\newcommand{\mcD}{\mathcal{ D}}
\newcommand{\mcE}{\mathcal{ E}}
\newcommand{\mcF}{\mathcal{ F}}

\newcommand{\mcI}{\mathcal{ I}}
\newcommand{\mcJ}{\mathcal{ J}}

\newcommand{\mcU}{\mathcal{ U}}
\newcommand{\mcV}{\mathcal{ V}}
\newcommand{\mcW}{\mathcal{ W}}

\newcommand{\mbP}{\mathbb{ P}}

\newcommand{\msS}{\mathscr{S}}

\begin{document}
\author{David Milovich}
\address{Department of Mathematics, University of Wisconsin, Madison, WI, 53706}
\email{milovich@math.wisc.edu}
\thanks{Support provided by an NSF graduate fellowship}

\begin{abstract}
Motivated by a question of Isbell, we show that $\Diamond$ implies there is a non\nbd-P\nbd-point $\mcU\in\beta\omega\setminus\omega$ such that neither $\la\mcU,\supseteq\ra$ nor $\la\mcU,\supseteq^*\ra$ is Tukey equivalent to $\la[\cardc]^{<\omega},\subseteq\ra$.  We also show that $\la\mcU,\supseteq^*\ra\equiv_T\la[\cardc]^{<\kappa},\subseteq\ra$ for some $\mcU\in\beta\omega\setminus\omega$, assuming $\cf(\kappa)=\kappa\leq\cardp=\cardc$.  We also prove two negative ZFC results about the possible Tukey classes of ultrafilters on $\omega$.
\end{abstract}

\title{Tukey classes of ultrafilters on $\omega$}

\date{July 25, 2008}

\maketitle

\section{Tukey classes}

\begin{definition} A quasiorder is a set with a transitive reflexive relation (denoted $\leq$ by default).  A quasiorder $Q$ is a $\kappa$\nbd-directed set if every subset of size less than $\kappa$ has an upper bound.  We abbreviate ``$\omega$\nbd-directed'' with ``directed.''
\end{definition}

\begin{definition} The product $P\times Q$ of two quasiorders $P$ and $Q$ is defined by $\la p_0,q_0\ra\leq\la p_1,q_1\ra$ iff $p_0\leq p_1$ and $q_0\leq q_1$.
\end{definition}

\begin{definition} A subset $C$ of a quasiorder $Q$ is cofinal if for all $q\in Q$ there exists $c\in C$ such that $q\leq c$.  The cofinality of $Q$ (written $\cf(Q)$), is defined as follows.
\begin{equation*}
\cf(Q)=\min\{\card{C}:C\text{ cofinal in }Q\}
\end{equation*}
\end{definition}

\begin{definition}[Tukey~\cite{tukey}] Given directed sets $P$ and $Q$ and a map $f\colon P\rightarrow Q$, we say $f$ is a Tukey map, writing $f\colon P\leq_T Q$, if the $f$\nbd-image of every unbounded subset of $P$ is unbdounded in $Q$.  We say $P$ is Tukey reducible to $Q$, writing $P\leq_T Q$, if there is a Tukey map from $P$ to $Q$.  If $P\leq_T Q\leq_T P$, then we say $P$ and $Q$ are Tukey equivalent and write $P\equiv_T Q$.
\end{definition}

\begin{proposition}[Tukey~\cite{tukey}] A map $f\colon P\rightarrow Q$ is Tukey if and only the $f$\nbd-preimage of every bounded subset of $Q$ is bounded in $P$.  Moreover, $P\leq_T Q$ if and only if there is a map $g\colon Q\rightarrow P$ such that the image of every cofinal subset of $Q$ is cofinal in $P$.  
\end{proposition}

\begin{theorem}[Tukey~\cite{tukey}]
$P\equiv_T Q$ if and only if $P$ and $Q$ order embed as cofinal subsets of a common third directed set.  Moreover, if $P\cap Q=\emptyset$, then we may assume the order embeddings are identity maps onto a quasiordering of $P\cup Q$.
\end{theorem}

The following is a list of basic facts about Tukey reducibility.
\begin{itemize}
\item $P\leq_T Q\Rightarrow \cf(P)\leq\cf(Q)$.
\item For all ordinals $\alpha,\beta$, we have $\alpha\leq_T\beta\Leftrightarrow\cf(\alpha)=\cf(\beta)$.
\item $P\leq_T P\times Q$.
\item $P\leq_T R\geq_T Q\Rightarrow P\times Q\leq_T R$.
\item $P\times P\equiv_T P$.
\item $P\leq_T\la[\cf(P)]^{<\omega},\subseteq\ra$.
\item For all infinite sets $A,B$, we have $\la[A]^{<\omega},\subseteq\ra\leq_T\la[B]^{<\omega},\subseteq\ra\Leftrightarrow\card{A}\leq\card{B}$.
\item Given finitely many ordinals $\alpha_0,\ldots,\alpha_{m-1},\beta_0,\ldots,\beta_{n-1}$, we have
\begin{equation*}
\prod_{i<m}\alpha_i\leq_T\prod_{i<n}\beta_i\Leftrightarrow\{\cf(\alpha_i):i<m\}\subseteq\{\cf(\beta_i):i<n\}.
\end{equation*}
\item Every countable directed set is Tukey equivalent to $1$ or $\omega$.
\end{itemize}

\begin{theorem}[Isbell~\cite{isbell}] No two of $1$, $\omega$, $\omega_1$, $\omega\times\omega_1$, and $\la[\omega_1]^{<\omega},\subseteq\ra$ are Tukey equivalent.  
\end{theorem}

Isbell~\cite{isbell} asked if these five Tukey classes encompass all directed sets of size $\omega_1$.  In \cite{isbell2}, he answered ``no'' assuming CH.  In particular, $\omega^\omega$, ordered by domination, is not Tukey equivalent to any of the above five orders.  Devlin, Stepr\=ans, and Watson~\cite{devlinsw} showed that $\Diamond$ implies there are $2^{\omega_1}$\nbd-many pairwise Tukey inequivalent directed sets of size $\omega_1$.  Todor\v cevi\`c~\cite{todorcevic} weakened the hypothesis of $\Diamond$ to CH and also showed that PFA implies that $1$, $\omega$, $\omega_1$, $\omega\times\omega_1$, and $\la[\omega_1]^{<\omega},\subseteq\ra$ represent the only Tukey classes of directed sets of size $\omega_1$.

\section{Tukey reducibility and topology}

Tukey reducibility has primarily been connected to topology by the concept of subnet: we say $\la x_i\ra_{i\in I}$ is a subnet of $\la y_j\ra_{j\in J}$ if there exists $f\colon I\rightarrow J$ such that the image of every cofinal subset of $I$ is cofinal in $J$, and $x_i=y_{f(i)}$ for all $i\in I$.  In contrast, our results are about classifying points in certain spaces by the Tukey classes of their local bases ordered by reverse inclusion.  The following theorem, which is of independent interest, implies that the Tukey class of a local base at a point in a space is a topological invariant.

\begin{theorem} Suppose $X$ and $Y$ are spaces, $p\in X$, $q\in Y$, $\mcA$ is a local base at $p$ in $X$, $\mcB$ is a local base at $q$ in $Y$, $f\colon X\rightarrow Y$ is continuous and open (or just continuous at $p$ and open at $p$), and $f(p)=q$.  Then $\la\mcB,\supseteq\ra\leq_T\la\mcA,\supseteq\ra$.
\end{theorem}
\begin{proof} Choose $H\colon\mcA\rightarrow\mcB$ such that $H(U)\subseteq f[U]$ for all $U\in\mcA$.  (Here we use that $f$ is open.)  Suppose $\mcC\subseteq\mcA$ is cofinal.  For any $U\in\mcB$, we may choose $V\in\mcA$ such that $f[V]\subseteq U$ by continuity of $f$.  Then choose $W\in\mcC$ such that $W\subseteq V$.  Hence, $H(W)\subseteq f[W]\subseteq f[V]\subseteq U$.  Thus, $H[\mcC]$ is cofinal.
\end{proof}

\begin{corollary}
In the above theorem, if $f$ is a homeomorphism, then every local base at $p$ is Tukey\nbd-equivalent to every local base at $q$.
\end{corollary}

\begin{example}
Consider the ordered space $X=\omega_1+1+\omega^{\text{op}}$.  It has a point $p$ that is the limit of an ascending $\omega_1$\nbd-sequence and a descending $\omega$-sequence.  Every local base at $p$, ordered by $\supseteq$, is Tukey equivalent to $\omega\times\omega_1$.

Next, consider $D_{\omega_1}\cup\{\infty\}$, the one-point compactification of the $\omega_1$\nbd-sized discrete space.  Glue $X$ and $D_{\omega_1}\cup\{\infty\}$ together into a new space $Y$ by a quotient map that identifies $p$ and $\infty$.  In $Y$, every local base at $p$, ordered by $\supseteq$, is Tukey equivalent to $\la [\omega_1]^{<\omega},\subseteq\ra$, which is not Tukey equivalent to $\omega\times\omega_1$.  

Thus, we can distinguish $p$ in $X$ from $p$ in $Y$ by their associated Tukey classes, even though other topological properties, such as character and $\pi$\nbd-character, have not changed.  Moreover, since $\omega\times\nolinebreak\omega_1<_T[\omega_1]^{<\omega}$, we may conclude there is no continuous open map from $X$ to $Y$ that sends $p$ to $p$.
\end{example}

\section{Ultrafilters}

\begin{definition}
Let $\omega^*$ denote the space $\beta\omega\setminus\omega$ of nonprincipal ultrafilters on $\omega$.
\end{definition}

By Stone duality, every ultrafilter $\mcU$ on $\omega$ is such that $\mcU$ ordered by containment, $\supseteq$, is Tukey-equivalent to every local base of $\mcU$ in $\beta\omega$.  Likewise, $\mcU$ ordered by almost containment, $\supseteq^*$, is Tukey equivalent to every local base of $\mcU$ in $\omega^*$.  
Therefore, let us now restrict our attention to the Tukey classes of nonprincipal ultrafilters on $\omega$, ordered by $\supseteq$ or $\supseteq^*$.  Note that the identity map on a $\mcU\in\omega^*$ is a Tukey map from $\la\mcU,\supseteq^*\ra$ to $\la\mcU,\supseteq\ra$.  Moreover, since $\la[\cardc]^{<\omega},\subseteq\ra$ is Tukey\nbd-maximal among the directed sets of cofinality at most $\cardc$, if $\la\mcU,\supseteq^*\ra\equiv_T\la[\cardc]^{<\omega},\subseteq\ra$, then $\la\mcU,\supseteq\ra\equiv_T\la[\cardc]^{<\omega},\subseteq\ra$.

Isbell~\cite{isbell}, using an independent family of sets, showed that there is always some $\mcU\in\omega^*$ such that $\la\mcU,\supseteq\ra\equiv_T\la[\cardc]^{<\omega},\subseteq\ra$.  Moreover, his proof also implicitly shows that $\la\mcU,\supseteq^*\ra\equiv_T\la[\cardc]^{<\omega},\subseteq\ra$.

\begin{definition}
We say $\mcI\subseteq[\omega]^\omega$ is independent if for all disjoint $\sigma,\tau\in[\mcI]^{<\omega}$ we have $\bigcap\sigma\not\subseteq^*\bigcup\tau$.
\end{definition}

\begin{lemma}[Hausdorff~\cite{hausdorff}] \label{LEMindepset}
There exists an independent $\mcI\in[[\omega]^\omega]^\cardc$.
\end{lemma}

\begin{theorem}[Isbell~\cite{isbell}]
There exists $\mcU\in\omega^*$ such that $\la\mcU,\supseteq^*\ra\equiv_T\la[\cardc]^{<\omega},\subseteq\ra$.
\end{theorem}
\begin{proof}
It suffices to show that there exists $f\colon\la[\cardc]^{<\omega},\subseteq\ra\leq_T\la\mcU,\supseteq^*\ra$.
Let $\mcI\in[[\omega]^\omega]^\cardc$ be independent.  Let $\mcF$ be the filter generated by $\mcI$.  Let $\mcJ$ be the ideal generated by the set of pseudointersections of infinite subsets of $\mcI$.  Extend $\mcF$ to an ultrafilter $\mcU$ disjoint from $\mcJ$.  Define $f\colon[\cardc]^{<\omega}\rightarrow\mcU$ by $\sigma\mapsto\bigcap_{\alpha\in\sigma}I_\alpha$.  Then $f$ is Tukey as desired.
\end{proof}

\begin{definition}
Given $\mcU\in\omega^*$, we say $\mcU$ is a $P_\kappa$\nbd-point if $\la\mcU,\supseteq^*\ra$ is $\kappa$\nbd-directed.  We call $P_{\omega_1}$\nbd-points P\nbd-points.
\end{definition}

There are also known constructions of various $\mcU\in\omega^*$ that satisfy $\la\mcU,\supseteq^*\ra\equiv_T\la[\cardc]^{<\omega},\subseteq\ra$ and some additional property.  See, for example, Dow and Zhou~\cite{dow}.  Also, Kunen~\cite{kunenok} proved that there exists a non\nbd-P\nbd-point $\mcU\in\omega^*$ such that $\mcU$ is $\cardc$\nbd-OK, and the next proposition shows that such a point must satisfy $\la\mcU,\supseteq^*\ra\equiv_T\la[\cardc]^{<\omega},\subseteq\ra$.

\begin{definition}[Kunen~\cite{kunenok}] \label{DEFkappaok}
We say $\mcU\in\omega^*$ is $\kappa$\nbd-OK if for every $\la A_n\ra_{n<\omega}\in\mcU^\omega$ there exists $\la B_\alpha\ra_{\alpha<\kappa}\in\mcU^{\kappa}$ such that for all nonempty $\sigma\in[\kappa]^{<\omega}$ we have $\bigcap_{\alpha\in\sigma}B_\alpha\subseteq^*A_{\card{\sigma}}$.
\end{definition}

\begin{proposition}
If $\mcU$ is a $\cardc$\nbd-OK non\nbd-P\nbd-point in $\omega^*$, then $\la\mcU,\supseteq^*\ra\equiv_T\la[\cardc]^{<\omega},\subseteq\ra$.
\end{proposition}
\begin{proof}
It suffices to show that there exists $f\colon\la[\cardc]^{<\omega},\subseteq\ra\leq_T\la\mcU,\supseteq^*\ra$.  Choose $\la A_n\ra_{n<\omega}\in\mcU^\omega$ such that $\{A_n:n<\omega\}$ has no pseudointersection in $\mcU$.  Then choose $\la B_\alpha\ra_{\alpha<\cardc}\in\mcU^{\cardc}$ as in Definition~\ref{DEFkappaok}.  Define $f\colon[\cardc]^{<\omega}\rightarrow \mcU$ by $\sigma\mapsto\bigcap_{\alpha\in\sigma}B_\alpha$.  Then every infinite subset of $[\cardc]^{<\omega}$ has unbounded $f$\nbd-image; hence, $f$ is Tukey as desired.
\end{proof}

\begin{definition}
Let $\cardu$ denote the least $\kappa$ such that there exists $\mcU\in\omega^*$ such that $\cf(\la\mcU,\supseteq^*\ra)=\kappa$.  Note that $\cf(\la\mcU,\supseteq\ra)=\cf(\la\mcU,\supseteq^*\ra)$ always holds.
\end{definition}

Isbell~\cite{isbell} asked if every $\mcU\in\omega^*$ satisfies $\la\mcU,\supseteq\ra\equiv_T\la[\cardc]^{<\omega},\subseteq\ra$.  It is now well\nbd-known that it is consistent with $\neg\text{CH}$ that $\cardu<\cardc$, which implies the existence of $\mcU\in\omega^*$ such that $\la\mcU,\supseteq\ra\leq_T\la[\cardu]^{<\omega},\subseteq\ra<_T\la[\cardc]^{<\omega},\subseteq\ra$.  To keep Isbell's question interesting, we must restrict our attention to models of $\cardu=\cardc$.

\begin{definition}
We say $\mcA\subseteq\powset{\omega}$ has the \emph{strong finite intersection property} (SFIP) if $\card{\bigcap\sigma}=\omega$ for all $\sigma\in[\mcA]^{<\omega}$.  Let $\cardp$ denote the least $\kappa$ for which some $\mcA\in[[\omega]^\omega]^\kappa$ has the SFIP but does not have a nontrivial pseudointersection. 
\end{definition}

It easily follows that $\cardp\leq\cardu$.  Moreover, by Bell's Theorem~\cite{bell}, $\cardp$ is the least $\kappa$ for which there exists a $\sigma$\nbd-centered poset $\mbP$ and a family $\mcD$ of $\kappa$\nbd-many dense subsets of $\mbP$ such that $\mbP$ does not have a $\mcD$\nbd-generic filter.  Hence, $\cardp=\cardc$ is equivalent to $\text{MA}_{\sigma\text{-centered}}$.  

\begin{definition}
Given cardinals $\kappa$ and $\lambda$, let $E^\kappa_\lambda$ denote $\{\alpha<\kappa:\cf(\alpha)=\lambda\}$.
\end{definition}

\begin{theorem} Assume $\Diamond(E^\cardc_\omega)$ and $\cardp=\cardc$.  Then there exists $\mcU\in\omega^*$ such that $\mcU$ is not a $P$\nbd-point and $\cardc<_T\la\mcU,\supseteq^*\ra\leq_T\la\mcU,\supseteq\ra<_T[\cardc]^{<\omega}$.
\end{theorem}
\begin{proof}
To simplify notation, we construct $\mcU$ as an ultrafilter on $\omega^2$.  Indeed, we construct $P_\cardc$\nbd-points $\mcV,\mcW_0,\mcW_1,\mcW_2,\ldots\in\omega^*$ and set $\mcU=\{E\subseteq\omega^2:\mcV\ni\{i:\mcW_i\ni\{j:\la i,j\ra\in E\}\}\}$.  This immediately implies that $\{(\omega\setminus n)\times\omega:n<\omega\}$ is a countable subset of $\mcU$ with no pseudointersection in $\mcU$; whence, $\mcU$ is not a $P$\nbd-point.  Our construction proceeds in $\cardc$ stages such that, for each $n<\omega$, the sequences $\la\mcV_\alpha\ra_{\alpha<\cardc}$ and $\la\mcW_{n,\alpha}\ra_{\alpha<\cardc}$ are continuous increasing chains of filters such that $\mcV=\bigcup_{\alpha<\cardc}\mcV_\alpha$ and $\mcW_n=\bigcup_{\alpha<\cardc}\mcW_{n,\alpha}$.  Set $\mcU_{\alpha}=\{E\subseteq\omega^2:\mcV_{\alpha}\ni\{i:\mcW_{i,\alpha}\ni\{j:\la i,j\ra\in E\}\}\}$ for all $\alpha<\cardc$.

Let $\la\Xi_\alpha\ra_{\alpha\in E^\cardc_\omega}$ be a $\Diamond$\nbd-sequence.  Let $\zeta\colon\cardc\leftrightarrow[\omega]^\omega$ and $\eta\colon\cardc\leftrightarrow\bigl[\omega^2\bigr]^\omega$.  Set $\mcV_0=\mcW_{n,0}=\{\omega\setminus\sigma:\sigma\in[\omega]^{<\omega}\}$ for all $n<\omega$.  Suppose $\alpha<\cardc$ and we've constructed $\la\mcV_\beta\ra_{\beta<\alpha}$ and $\la\mcW_{n,\beta}\ra_{\la n,\beta\ra\in\omega\times\alpha}$ such that, for all $\beta<\alpha$ and $n<\omega$, $\mcV_\beta$ and $\mcW_{n,\beta}$ are filters on $\omega$; if $\cf(\beta)\not=\omega$ and $\beta+1<\alpha$, then further suppose that $\mcV_\beta$ and $\mcW_{n,\beta}$ have pseudointersections in $\mcV_{\beta+1}$ and $\mcW_{n,\beta+1}$, respectively.  If $\alpha$ is a limit ordinal, then set $\mcV_\alpha=\bigcup_{\beta<\alpha}\mcV_\beta$ and $\mcW_{n,\alpha}=\bigcup_{\beta<\alpha}\mcW_{n,\beta}$ for each $n<\omega$.  If $\alpha$ is the successor of an ordinal with cofinality other than $\omega$, then we use stage $\alpha$ as follows to help our filters become ultrafilters that are $P_\cardc$\nbd-points.  Choose the least $\beta<\cardc$ such that $\zeta(\beta),\omega\setminus\zeta(\beta)\not\in\mcV_{\alpha-1}$.  Choose $E\in\{\zeta(\beta),\ \omega\setminus\zeta(\beta)\}$ such that $\{E\}\cup\mcV_{\alpha-1}$ has the SFIP and let $\mcV_\alpha$ be a filter generated by $\mcV_{\alpha-1}$ and a pseudointersection of $\{E\}\cup\mcV_{\alpha-1}$.  Likewise, for each $n<\omega$, choose the least $\beta<\cardc$ such that $\zeta(\beta),\omega\setminus\zeta(\beta)\not\in\mcW_{n,\alpha-1}$.  Choose $E\in\{\zeta(\beta),\ \omega\setminus\zeta(\beta)\}$ such that $\{E\}\cup\mcW_{n,\alpha-1}$ has the SFIP and let $\mcW_{n,\alpha}$ be a filter generated by $\mcW_{n,\alpha-1}$ and a pseudointersection of $\{E\}\cup\mcW_{n,\alpha-1}$.

Finally, suppose $\alpha$ is the successor of an ordinal with cofinality $\omega$.  Then we use stage $\alpha$ to kill a potential witness to $\la\mcU,\supseteq\ra\equiv_T[\cardc]^{<\omega}$.  Choose, if it exists, the least $\beta<\cardc$ for which $\eta(\beta)$ is contained in the intersection of an infinite subset of $\eta[\Xi_\alpha]$ and $\{\eta(\beta)\}\cup\mcU_{\alpha-1}$ has the SFIP.  Let $\mcV_\alpha$ be the filter generated by $\{F\}\cup\mcV_{\alpha-1}$ where $F=\{i:\mcW_{i,\alpha-1}\not\ni\omega\setminus\{j:\la i,j\ra\in\eta(\beta)\}\}$; for each $i\in F$, let $\mcW_{i,\alpha}$ be the filter generated by $\{\{j:\la i,j\ra\in\eta(\beta)\}\}\cup\mcW_{i,\alpha-1}$; for each $i\in\omega\setminus F$, set $\mcW_{i,\alpha}=\mcW_{i,\alpha-1}$.  Note that this implies $\eta(\beta)\in\mcU_\alpha$.  If no such $\beta$ exists, then set $\mcV_\alpha=\mcV_{\alpha-1}$ and $\mcW_{n,\alpha}=\mcW_{n,\alpha-1}$ for all $n<\omega$.  This completes the construction.

Clearly, $\cardc\leq_T\la\mcV,\supseteq^*\ra\leq_T\la\mcU,\supseteq^*\ra$.  Since $\mcU$ is not a $P$\nbd-point, $\cardc\not\equiv_T\la\mcU,\supseteq^*\ra$.
Therefore, it remains only to show that  $\la\mcU,\supseteq\ra\not\equiv_T[\cardc]^{<\omega}$.  Suppose $\mcA\in[\mcU]^\cardc$.  Then it suffices to show that the intersection of an infinite subset of $\mcA$ is in $\mcU$.  By $\Diamond(E^\cardc_\omega)$, there exists $M\elemsub H_{\cardc^+}$ such that $\card{M}=\omega$ and $M\supseteq\{\mcA,\la\mcV_\alpha\ra_{\alpha<\cardc},\la\mcW_{n,\alpha}\ra_{\la n,\alpha\ra\in\omega\times\cardc}\}$ and $\eta[\Xi_\delta]=\mcA\cap M$ where $\delta=\sup(\cardc\cap M)$.  Hence, it suffices to show that the intersection $E$ of some infinite subset of $\mcA\cap M$ is such that $\{E\}\cup\mcU_{\delta}$ has the SFIP.  

Let $\{V_n:n<\omega\}$ generate of the filter $\mcV_\delta$; for each $i<\omega$, let $\{W_{i,j}:j<\omega\}$ generate the filter $\mcW_{i,\delta}$.  Set $\mcB_0=\mcA$.  Suppose $k<\omega$ and, for all $l<k$, we have $A_l\in\mcB_{l+1}\in[\mcB_l]^{\cardc}$ and $n_l<\omega$ and $\mcW_{n_l}\ni\bigl\{j:\la n_l,j\ra\in B\cap\bigcap_{h<l}A_h\bigr\}$ for all $B\in\mcB_{l+1}$.
Since $\cf(\cardc)>\omega$, there exist $\mcB_{k+1}\in[\mcB_k]^{\cardc}$ and $n_k\in\bigcap_{h<k}(V_h\setminus\{n_h\})$ and $\sigma_k\colon\{n_l:l<k\}\rightarrow\omega$ such that, for all $l<k$ and $B\in\mcB_{k+1}$, we have $\mcW_{n_k}\ni\bigl\{j:\la n_k,j\ra\in B\cap\bigcap_{h<k}A_h\bigr\}$ and $\sigma_k(n_l)\in\bigcap_{h<k}W_{n_l,h}$ and $\sigma_k\subseteq B\cap\bigcap_{h<k}A_h$.  Choose any $A_k\in\mcB_{k+1}\setminus\{A_h:h<k\}$.  By induction, we can repeat the above for all $k<\omega$.  Moreover, we may carry out any finite initial segment of the construction in $M$.  Hence, \wma\ $\{A_i:i<\omega\}\subseteq M$.  Finally, $\bigcup_{i<\omega}\sigma_i\subseteq\bigcap_{i<\omega}A_i$ and $\{\bigcup_{i<\omega}\sigma_i\}\cup\mcU_\delta$ has the SFIP.
\end{proof}

Note that $\Diamond(E^\cardc_\omega)$ is equivalent to $\Diamond$ under CH.  Furthermore, a recent result of Shelah~\cite{shelah} is that if $\kappa$ is an uncountable cardinal and $2^\kappa=\kappa^+$, then $\Diamond(S)$ holds for every stationary $S$ disjoint from $E^{\kappa^+}_{\cf(\kappa)}$.  Hence, we could drop the hypothesis $\Diamond(E^\cardc_\omega)$ under the assumption that $\cardc=\kappa^+$ for some cardinal $\kappa$ of uncountable cofinality.  (We'd have $2^\kappa=\kappa^+$ because $\cardc^{<\cardp}=\cardc$. (See Martin and Solovay~\cite{martinsolovay}.)

It is worth noting another relationship between the Tukey classes arising from ultrafilters ordered by $\supseteq^*$ and those ordered by $\supseteq$.

\begin{proposition}\label{PROnonptukmap}
Suppose $\mcU$ is a non\nbd-P\nbd-point in $\omega^*$.  Then there exists $\mcV\in\omega^*$ such that $\la\mcV,\supseteq\ra\leq_T\la\mcU,\supseteq^*\ra$.
\end{proposition}
\begin{proof}
Choose $\la x_n\ra_{n<\omega}\in\mcU^\omega$ such that $x_n\supseteq x_{n+1}\not\supseteq^*x_n$ for all $n<\omega$, that $\bigcap_{n<\omega}x_n=\emptyset$, and that $\{x_n:n<\omega\}$ has no pseudointersection in $\mcU$.  For each $n<\omega$, set $y_n=x_n\setminus x_{n+1}$.  Set $\mcV=\left\{E\subseteq\omega: \bigcup_{n\in E}y_n\in\mcU\right\}$.  Then $\mcV\in\omega^*$ and the map from $\la\mcV,\supseteq\ra$ to $\la\mcU,\supseteq^*\ra$ defined by $E\mapsto\bigcup_{n\in E}y_n$ is Tukey.
\end{proof}

Next, we have a pair of negative ZFC results.

\begin{theorem}\label{THMMAomegaQ}
Let $Q$ be a directed set that is a countable union of $\omega_1$\nbd-directed sets.  Then $\la\mcU,\supseteq^*\ra\not\equiv_T\omega\times Q$ for all $\mcU\in\omega^*$.
\end{theorem}
\begin{proof}
Seeking a contradiction, suppose $\mcU\in\omega^*$ and $\la\mcU,\supseteq^*\ra\equiv_T\omega\times Q$.  Then there is a quasiordering $\sqsubseteq$ on $\mcU\cup(\omega\times Q)$ such that $\la\mcU,\supseteq^*\ra$ and $\la\omega\times Q,\leq_{\omega\times Q}\ra$ are cofinal suborders.  Let $Q=\bigcup_{n<\omega}Q_n$ where $Q_n$ is $\omega_1$\nbd-directed for all $n<\omega$.  Fix $p\in Q$.  Fix $\eta\in\omega^\omega$ such that $\inv{\eta}\{n\}$ is unbounded and $\eta(4n)=\eta(4n+1)=\eta(4n+2)=\eta(4n+3)$ for all $n<\omega$.  For all $n<\omega$ and $q\in Q$, choose $x_{n,q}\in\mcU$ such that $\la n,q\ra\sqsubseteq x_{n,q}$.  \Wma\ that $x_{i,p}\sqsubseteq x_{j,q}$ for all $i\leq j<\omega$ and $q\in Q$.

Construct $\zeta\in\omega^\omega$ as follows.  Suppose we are given $n<\omega$ and $\zeta\restrict n$.  Then, for all $q\in Q$, the set $\{x_{\zeta(m),q}:m<n\}$ has a $\sqsubseteq$\nbd-upper bound $\la k,r\ra$ for some $k<\omega$ and $r\in Q$.  Since $Q_{\eta(n)}$ is $\omega_1$\nbd-directed, every countable partition of $Q_{\eta(n)}$ includes a cofinal subset.  Hence, there exist $k<\omega$ and a cofinal subset $S_n$ of $Q_{\eta(n)}$ such that for all $q\in S_n$ there exists $r\in Q$ such that $\{x_{\zeta(m),q}:m<n\}\sqsubseteq\la k,r\ra$.  \Wma\ $k>\zeta(m)$ for all $m<n$. Set $\zeta(n)=k$.

Since $\omega^*$ is an F\nbd-space (or, more directly, by an easy diagonalization argument), there exists $z\subseteq\omega$ such that $x_{\zeta(4n),p}\setminus x_{\zeta(4n+2),p}\subseteq^*z$ and $x_{\zeta(4n+2),p}\setminus x_{\zeta(4n+4),p}\subseteq^*\omega\setminus z$ for all $n<\omega$.  Suppose $z\in\mcU$.  Then there exist $m<\omega$ and $\la l,r\ra\in\omega\times Q_m$ such that $\la l,r\ra\sqsupseteq z$.  Choose $n<\omega$ such that $\eta(4n+3)=m$ and $\zeta(4n+2)\geq l$.  Then choose $q\in S_{4n+3}$ such that $q\geq r$.  Then $\la\zeta(4n+2),q\ra\sqsupseteq z$. Hence, $x_{\zeta(4n+2),q}\sqsupseteq z\cap x_{\zeta(4n+2),p}\sqsupseteq x_{\zeta(4n+4),p}\sqsupseteq\la\zeta(4n+4),p\ra$.  Hence, $\la\zeta(4n+4),p\ra\sqsubseteq x_{\zeta(4n+2),q}\sqsubseteq\la\zeta(4n+3),s\ra$ for some $s\in Q$, which is absurd because $\zeta$ is strictly increasing.  By symmetry, we can also derive an absurdity from $\omega\setminus z\in\mcU$.  Thus, $\mcU$ is not an ultrafilter on $\omega$, which yields our desired contradiction.
\end{proof}

The above result is optimal in the following sense.  It is not hard to show that, for a fixed regular uncountable $\kappa$, a construction of Brendle and Shelah~\cite{brendle} can be trivially modified to yield of a model of ZFC in which some $\mcU\in\omega^*$ satisfies $\la\mcU,\supseteq^*\ra\equiv_T\kappa\times\lambda$ for each $\lambda$ in an arbitrary set of regular cardinals exceeding $\kappa$.

\begin{definition}
A quasiorder $Q$ is said to be $\kappa$\nbd-like if every bounded subset of $Q$ has size less than $\kappa$.
\end{definition}

\begin{lemma}\label{LEMcfoplike}
Given a quasiorder $Q$ with an unbounded cofinal subset $C$, there exists a cofinal subset $A$ of $C$ such that $A$ is $\card{C}$\nbd-like.
\end{lemma}
\begin{proof}
Let $\la c_\alpha\ra_{\alpha<\card{C}}\colon\card{C}\leftrightarrow C$.  For each $\alpha<\card{C}$, let $a_\alpha=c_\beta$ where $\beta$ is the least $\gamma<\card{C}$ such that $c_\gamma$ has no upper bound in $\{a_\delta:\delta<\alpha\}$, provided such a $\gamma$ exists.  If no such $\gamma$ exists, then $\alpha>0$, so we may set $a_\alpha=a_0$.  Then $A=\{a_\alpha:\alpha<\card{C}\}$ is as desired.
\end{proof}

\begin{theorem}
Suppose $Q$ is a directed set that is a countable union of $\omega_1$\nbd-directed sets.  Then $\la\mcU,\supseteq\ra\not\leq_T Q$ for all $\mcU\in\omega^*$ satisfying $\cf(\cf(\la\mcU,\supseteq\ra))>\omega$.
\end{theorem}
\begin{proof}
Seeking a contradiction, suppose $\mcU\in\omega^*$ and $\cf(\cf(\la\mcU,\supseteq\ra))>\omega$ and $f:\la\mcU,\supseteq\ra\leq_T Q$.  By Lemma~\ref{LEMcfoplike}, $\mcU$ has a cofinal subset $\mcA$ that is $\cf(\la\mcU,\supseteq\ra)$\nbd-like.  Since $\mcA$ is cofinal, $f\restrict\mcA$ is a Tukey map and $\card{\mcA}=\cf(\la\mcU,\supseteq\ra)$.  Let $Q=\bigcup_{n<\omega}Q_n$ where $Q_n$ is $\omega_1$\nbd-directed for all $n<\omega$.  Since $\cf(\card{\mcA})>\omega$, there exist $n<\omega$ and $\mcB\in[\mcA]^{\card{\mcA}}$ such that $f[\mcB]\subseteq Q_n$.  Since $\mcA$ is $\card{\mcA}$\nbd-like, $\mcB$ is unbounded.  Set $I=\omega\setminus\bigcap\mcB$.  For each $i\in I$, choose $B_i\in\mcB$ such that $i\not\in B_i$.  Then $\bigcap_{i\in I}B_i=\bigcap\mcB$; hence, $\{B_i:i\in I\}$ is unbounded.  But $\{f(B_i):i\in I\}$ is a countable subset of $Q_n$, and therefore bounded.  This contradicts our assumption that $f$ is Tukey.
\end{proof}

Our next theorem is a positive consistency result.  Its proof uses Solovay's Lemma~\cite{martinsolovay}, which we now state in terms of $\cardp$.

\begin{lemma}
If $\mcA,\mcB\in[[\omega]^\omega]^{<\cardp}$ and $\card{a\cap\bigcap\sigma}=\omega$ for all $a\in\mcA$ and $\sigma\in[\mcB]^{<\omega}$, then $\mcB$ has a pseudointersection $b$ such that $\card{a\cap b}=\omega$ for all $a\in\mcA$.
\end{lemma}

\begin{theorem}  Assume $\cardp=\cardc$.  Let $\omega\leq\cf(\kappa)=\kappa\leq\cardc$.  Then there exists $\mcU\in\omega^*$ such that $\la\mcU,\supseteq^*\ra\equiv_T\la[\cardc]^{<\kappa},\subseteq\ra$.
\end{theorem}
\begin{proof}
Given a set $E$, let $I(E)$ denote the set of injections from $\kappa$ to $E$.  Given $\mcE\subseteq\powset{\omega}$, let $\Phi(\mcE)$ denote the set of $\la\rho,\Gamma\ra\in[\mcE]^{<\omega}\times I(\mcE)^{<\omega}$ satisfying $\bigcap\rho\subseteq^*\bigcup_{f\in\ran\Gamma}f(\gamma)$ for all $\gamma<\kappa$.  Let $\la\msS_\alpha\ra_{\alpha<\cardc}$ enumerate $[[\omega]^\omega]^{<\kappa}$.  Note that if $\card{\mcE}\geq\kappa$, then $\Phi(\mcE)=\emptyset$ implies that $\mcE$ has the SFIP and that $\la\mcE,\supseteq^*\ra$ is $\kappa$\nbd-like.

Let us construct a sequence $\la U_\alpha\ra_{\alpha<\cardc}$ in $[\omega]^\omega$ such that we have the following for all $\alpha\leq\cardc$, given the notation $\mcU_\beta=\{U_\gamma:\gamma<\beta\}$ for all $\beta\leq\cardc$.
\begin{enumerate}
\item\label{enumsplit}$\forall\beta<\alpha\ \,\forall\sigma,\tau\in[\mcU_\beta]^{<\omega}\ \, \bigcap\sigma\subseteq^*\bigcup\tau$ or $\bigcap\sigma\setminus\bigcup\tau\not\subseteq^*U_\beta$
\item\label{enumpkappa}$\forall\beta<\alpha\ \,\exists\sigma\in[\msS_\beta]^{<\omega}\ \,U_\beta\cap\bigcap\sigma=^*\emptyset$ or $\forall S\in\msS_\beta\ \,U_\beta\subseteq^*S$
\item\label{enumphi}$\Phi(\mcU_\alpha)=\emptyset$
\end{enumerate}

Clearly, (\ref{enumsplit}) and (\ref{enumpkappa}) will be preserved at limit stages of the construction.  Let us show that (\ref{enumphi}) will also be preserved.  Let $\omega\leq\cf(\eta)\leq\eta\leq\cardc$ and suppose (\ref{enumsplit}) and (\ref{enumphi}) hold for all $\alpha<\eta$.  Seeking a contradiction, suppose $\la\rho,\Gamma\ra\in\Phi(\mcU_\eta)$; \wma\ $\la\rho,\Gamma\ra$ is chosen so as to minimize $\dom\Gamma$.  By (\ref{enumsplit}), $\la U_\alpha\ra_{\alpha<\eta}$ is injective; let $\psi$ be its inverse.  Since $\Phi(\mcU_{\sup(\psi[\rho])})=\emptyset$, we have $\Gamma\not=\emptyset$.  By the pigeonhole principle, there exist $A\in[\kappa]^\kappa$ and $i\in\dom\Gamma$ such that for all $\gamma\in A$ we have $\psi(\Gamma(i)(\gamma))=\max_{j\in\dom\Gamma}\psi(\Gamma(j)(\gamma))$. By symmetry, \wma\ $i=\max(\dom\Gamma)$.  Since $\Phi(\mcU_{\sup(\psi[\rho])})=\emptyset$, we have $\card{A\cap\inv{\Gamma(i)}\sup(\psi[\rho])}<\kappa$; hence, \wma\ $A\cap\inv{\Gamma(i)}\sup(\psi[\rho])=\emptyset$.  By the definition of $\Phi(\mcU_\eta)$, we have $\bigcap\rho\setminus\bigcup_{j<i}\Gamma(j)(\gamma)\subseteq^*\Gamma(i)(\gamma)$ for all $\gamma\in A$.  Hence, by (\ref{enumsplit}), we have $\bigcap\rho\subseteq^*\bigcup_{j<i}\Gamma(j)(\gamma)$ for all $\gamma\in A$.  Choose $h\in I(A)$.  Then $\la\rho,\la\Gamma(j)\circ h\ra_{j<i}\ra\in\Phi(\mcU_\eta)$, in contradiction with the minimality of $\dom\Gamma$.  Thus, (\ref{enumphi}) will be preserved at limit stages.

Given $\alpha<\cardc$ and $\la U_\beta\ra_{\beta<\alpha}$ satisfying (1)\nbd-(3), let us show that there always exists $U_\alpha\in[\omega]^\omega$ such that $\la U_\beta\ra_{\beta\leq\alpha}$ also satifies (1)\nbd-(3).  Let $g\in 2^\omega$ be sufficiently Cohen generic.  There are two cases to consider.  First, suppose that there exists $\sigma\in[\msS_\alpha]^{<\omega}$ such that $\Phi(\mcU_\alpha\cup\sigma)\not=\emptyset$.  Then there exists $\la\rho_2,\Gamma_2\ra\in\Phi(\mcU_\alpha\cup\{x_2\})$ where $x_2=\bigcap\sigma$.  For each $i<2$, set $x_i=\inv{g}\{i\}\setminus x_2$.  Seeking a contradiction, suppose there exists $\la\rho_i,\Gamma_i\ra\in\Phi(\mcU_\alpha\cup\{x_i\})$ for each $i<2$.  \Wma\ $\bigcup_{i<3}\bigcup\ran\Gamma_i\subseteq\mcU_\alpha$.  Let $\Lambda$ be a concatenation of $\{\Gamma_i:i<3\}$ and set $\tau=\mcU_\alpha\cap\bigcup_{i<3}\rho_i$.  Then, for all $\gamma<\kappa$, we have
\begin{equation*}
\bigcap\tau=\bigcup_{i<3}\left(x_i\cap\bigcap\tau\right)\subseteq\bigcup_{i<3}\bigcap\rho_i\subseteq^*\bigcup_{f\in\ran\Lambda}f(\gamma).
\end{equation*}
Hence, $\la\tau,\Lambda\ra\in\Phi(\mcU_\alpha)$, in contradiction with (\ref{enumphi}).  Therefore, we may choose $i<2$ such that $\Phi(\mcU_\alpha\cup\{x_i\})=\emptyset$.  Set $U_\alpha=x_i$, which is disjoint from $\bigcap\sigma$.  Then (\ref{enumpkappa}) and (\ref{enumphi}) are clearly satisfied for stage $\alpha+1$, and (\ref{enumsplit}) is also satisfied because of Cohen genericity.

Now suppose that $\Phi(\mcU_\alpha\cup\sigma)=\emptyset$ for all $\sigma\in[\msS_\alpha]^{<\omega}$.  For each $\rho\in[\mcU_\alpha]^{<\omega}$, $\sigma\in[\msS_\alpha]^{<\omega}$, and $\Gamma\in I(\mcU_\alpha)^{<\omega}$, choose $\gamma_{\rho,\sigma,\Gamma}<\kappa$ such that $\bigcap(\rho\cup\sigma)\not\subseteq^*\bigcup_{i\in\ran\Gamma}f(\delta)$ for all $\delta\in\kappa\setminus\gamma_{\rho,\sigma,\Gamma}$.  Set $\gamma_{\rho,\Gamma}=\sup\{\gamma_{\rho,\sigma,\Gamma}:\sigma\in[\msS_\alpha]^{<\omega}\}$; set $x_{\rho,\Gamma}=\bigcap\rho\setminus\bigcup_{f\in\ran\Gamma}f(\gamma_{\rho,\Gamma})$.  Then $x_{\rho,\Gamma}\cap\bigcap\sigma$ is infinite for all $\sigma\in[\msS_\alpha]^{<\omega}$.  By Solovay's Lemma, $\msS_\alpha$ has a pseudointersection $y$ such that $y\cap x_{\rho,\Gamma}$ is infinite for all $\rho\in[\mcU_\alpha]^{<\omega}$ and $\Gamma\in I(\mcU_\alpha)^{<\omega}$, for there are at most $\card{\mcU_\alpha}^{<\omega}$\nbd-many possible $x_{\rho,\Gamma}$.  Set $U_\alpha=y\cap\inv{g}\{0\}$.  Then (\ref{enumpkappa}) is clearly satisfied for stage $\alpha+1$.  Since $y\cap x_{\rho,\Gamma}\cap\bigcap\sigma$ is infinite, Cohen genericity implies $U_\alpha\cap x_{\rho,\Gamma}$ is infinite, for all $\rho$, $\sigma$, and $\Gamma$.  Hence, (\ref{enumphi}) is satisfied for stage $\alpha+1$; (\ref{enumsplit}) is also satisfied because of Cohen genericity.  This completes our construction of $\la U_\alpha\ra_{\alpha<\cardc}$.

Let $\mcU$ be the semifilter generated by $\mcU_\cardc$.  By (\ref{enumphi}), $\mcU_\cardc$ has the SFIP and $\mcU_\cardc$ is $\kappa$\nbd-like with respect to $\supseteq^*$.  Hence, by (\ref{enumpkappa}), $\mcU$ is a $P_\kappa$\nbd-point in $\omega^*$. Therefore, $f\colon\la\mcU,\supseteq^*\ra\leq_T\la[\cardc]^{<\kappa},\subseteq\ra$ for any injection $f$ of $\mcU$ into $[\cardc]^1$.  Choose $\zeta\colon[\cardc]^{<\kappa}\rightarrow\mcU$ such that $\zeta(\sigma)$ is a pseudointersection of $\{U_\alpha:\alpha\in\sigma\}$ for all $\sigma\in[\cardc]^{<\kappa}$.  Then $\zeta$ is Tukey because $\mcU_\cardc$ is $\kappa$\nbd-like.  Thus, $\mcU\leq_T[\cardc]^{<\kappa}\leq_T\mcU$.
\end{proof}

\section{Questions}

\begin{question}  Is it consistent that every $\mcU\in\omega^*$ satisfies $\la\mcU,\supseteq^*\ra\equiv_T\la[\cardc]^{<\omega}\ra$?
\end{question}

\begin{question}  Does CH (or even ZFC alone) imply there exists $\mcU\in\omega^*$ such that $\la\mcU,\supseteq\ra<_T\la[\cardc]^{<\omega}\ra$?
\end{question}

\begin{question} Does CH (or even ZFC alone) imply there exists a non\nbd-P\nbd-point $\mcU\in\omega^*$ such that $\la\mcU,\supseteq^*\ra<_T\la[\cardc]^{<\omega}\ra$? By Proposition~\ref{PROnonptukmap}, a positive answer to this question implies a positive answer to the previous question.
\end{question}

\begin{question} Does $\Diamond$ imply there are at least three Tukey classes represented by $\la\mcU,\supseteq^*\ra$ for some $\mcU\in\omega^*$?  Infinitely many Tukey classes?  As many as $2^{\omega_1}$?  What if we replace $\supseteq^*$ with $\supseteq$?
\end{question}

\begin{question}
Is it consistent with $\omega_1<\cardp$ that there exists $\mcU\in\omega^*$ such that $\la\mcU,\supseteq^*\ra\equiv_T\omega_1\times\cardc$?
\end{question}

\begin{question}
Does there consistently exist $\mcU\in\omega^*$ such that $$\cf(\cf(\la\mcU,\supseteq\ra)=\omega?$$
\end{question}

\begin{question}
Does there exist $\mcU\in\omega^*$ such that $\la\mcU,\supseteq\ra\equiv_T\omega^\omega$ where $\omega^\omega$ is ordered by domination?
\end{question}

\bibliographystyle{plain}

\end{document}